\documentclass[twoside,11pt]{amsart}
\usepackage[active]{srcltx} 
\usepackage[all]{xy}
\CompileMatrices
\usepackage{amsfonts}
\usepackage{amssymb}
\usepackage{amsthm}
\usepackage{amsmath}
\usepackage{ifthen}
\usepackage{url}
\usepackage{scalerel}
\usepackage{multicol}
\usepackage{multirow}
 \newtheorem{theorem}{Theorem}[section]
 \newtheorem{proposition}[theorem]{Proposition}
 \newtheorem{corollary}[theorem]{Corollary}
 \newtheorem{lemma}[theorem]{Lemma}
  \newtheorem{remark}[theorem]{Remark}

\theoremstyle{definition}
\newtheorem{definition}[theorem]{Definition}

\newtheorem{example}[theorem]{Example}

\DeclareMathOperator{\Gal}{Gal}            
\DeclareMathOperator{\Norm}{\textit{N}}    
\DeclareMathOperator{\Log}{Log}            

\DeclareMathOperator{\ab}{ab}
\DeclareMathOperator{\Pl}{Pl}
\DeclareMathOperator{\nr}{nr}

\DeclareMathOperator{\lr}{lr}
\DeclareMathOperator{\rank}{rank}

\DeclareMathOperator{\lc}{lc}

\DeclareMathOperator{\im}{im}

\newcommand{\Q}{\mathbb{Q}}
\newcommand{\F}{\mathbb{F}}
\newcommand{\Z}{\mathbb{Z}}
\newcommand{\N}{\mathbb{N}}

\newcommand{\Zl}{\mathbb{Z}_{\ell}}

\newcommand{\Ql}{\mathbb{Q}_{\ell}}
\newcommand{\Qp}{\mathbb{Q}_{p}}
\newcommand{\LL}{\Lambda}
\newcommand{\Kp}{K_{\p}}
\newcommand{\Lp}{L_{\mathfrak{P}}}

\newcommand{\p}{\mathfrak{p}}
\newcommand{\pp}{\mathfrak{P}}

\newcommand{\x}{\chi}

\newcommand{\Jc}{\mathcal{J}}
\newcommand{\Rc}{\mathcal{R}}
\newcommand{\Uc}{\mathcal{U}}
\newcommand{\Cc}{\mathcal{C}}

\newcommand{\Clog}[1]{\widetilde{C\ell}_{#1}}

\newcommand{\lval}[1]{\widetilde{v}_{#1}}
\newcommand{\mut}{\widetilde{\mu}}
\newcommand{\lat}{\widetilde{\lambda}}

\newcommand{\f}{\varphi}
\newcommand{\w}{\omega}

\newcommand{\isom}{\simeq}

\newcommand\scale[2]{\vstretch{#1}{\hstretch{#1}{#2}}}
\newcommand\ssp[2]{#1^{\scale{.6}{\scriptstyle {#2}}}}

\def\sectionnam{\@empty}
\def\subsectionnam{\@empty}

\begin{document}

\title[Iwasawa Logarithmic Invariants]{On the mu and lambda invariants of the logarithmic class group.}%

\author{Jos\'e-Ibrahim Villanueva-Guti\'errez}%
\address{Universit\'e de Bordeaux\\ Institut de Math\'ematiques\\
351, cours de la Lib\'eration \\F 33 405 TALENCE, France.} \email{jovillan@math.u-bordeaux.fr}
\urladdr{http://www.math.u-bordeaux.fr/\textasciitilde jovillan/}
\thanks{This paper was written under the financial support of CONACYT (The Mexican Council of Science and Technology) as part of the author's PhD project.}

\date{\today}%

\maketitle
\thispagestyle{empty}

\begin{abstract}
Let $\ell$ be a rational prime number. Assuming the Gross-Kuz'min conjecture along a $\Zl$-extension $K_{\infty}$ of a number field $K$, we show that there exist integers $\mut$, $\lat$ and $\widetilde{\nu}$ such that the exponent $\tilde{e}_{n}$ of the order $\ell^{\tilde{e}_{n}}$ of the logarithmic class group $\Clog{n}$ for the $n$-th layer $K_{n}$ of $K_{\infty}$ is given by $\tilde{e}_{n}=\widetilde{\mu}\ell^{n}+\widetilde{\lambda} n + \widetilde{\nu}$, for $n$ big enough. We show some relations between the classical invariants $\mu$ and $\lambda$, and their logarithmic counterparts $\mut$ and $\lat$ for some class of $\Zl$-extensions. Additionally, we provide numerical examples for the cyclotomic and the non-cyclotomic case. 
\end{abstract} 

\section{Introduction}

Let $K$ be a number field and $\ell$ a prime integer. The logarithmic class group $\Clog{K}$ of $K$ is the $\Zl$-module measuring the gap	 between the cyclotomic $\Zl$-extension $K^{\textrm{c}}$ of $K$ and the maximal abelian pro-$\ell$-extension $K^{\lc}$ of $K$ which splits completely over $K^{c}$, i.e. $\Clog{K}$ is isomorphic to the relative Galois group $\Gal(K^{\lc}/K^{c})$ (e.g. \cite{Jaulent86}, \cite[Thm. \& Conj. 2.3]{Jaulent94}). The logarithmic class group is conjectured to be finite for every number field, this is equivalent to the Gross-Kuz'min conjecture \cite[\S 2]{Jaulent17}; its finiteness is verified in the abelian case and some other cases (e.g. \cite{Jaulent16Naifs,Jaulent02}). As an arithmetic invariant of a number field the logarithmic class group has importance in its own. The logarithmic class group is effective by computational methods \cite{Diaz&Soriano99,Diaz&Jaulent&Pauli&Soriano05,Belabas&Jaulent16} and it is related to the wild kernels in $K$-theory  \cite{Jaulent&Soriano01}. In spite of its tight relation to the $\ell$-class group $C\ell_{K}$ of $K$, it behaves differently in several situations (see \S \ref{sec:numericalex}). 

We are interested in the study of the logarithmic class group along $\Zl$-extensions, in the spirit of Iwasawa's work for class groups \cite{Iwasawa59Gamma,Iwasawa73}. This has been done in the cyclotomic case by Jaulent in \cite{Jaulent86}; in this situation the logarithmic class group attached to the $n$-th layer of the cyclotomic \mbox{$\Zl$-extension} has a nice interpretation as a quotient of Kuz'min-Tate modules \cite{Jaulent17}. This interpretation is no longer valid for non-cyclotomic \mbox{$\Zl$-extensions} whenever the base field admits such an extension. 

Assuming the Gross-Kuz'min conjecture, the main result of this work provides a method for the non-cyclotomic case, and therefore completes the study of the logarithmic class group along $\Zl$-extensions:

\begin{theorem}\label{thm:IwasClog}
Let $K_{\infty}/K$ be a $\Zl$-extension. Assume the Gross-Kuz'min conjecture is valid along the \mbox{$\Zl$-extension} $K_{\infty}$. Let $\Clog{n}$ be the logarithmic class group of $K_{n}$ and let $\ell^{\widetilde{e}_{n}}$ be its order. There exist integers $\lat,\mut \geq 0$ and $\widetilde{\nu}$ such that 
\begin{equation*}
\widetilde{e}_{n}=\mut \ell^{n}+\lat n +\widetilde{\nu},\	\	\	\textrm{ for } n \textrm{ big enough.}
\end{equation*}
\end{theorem}

This article is organized as follows. In Section 2 we recall briefly the basic definitions and properties of the logarithmic arithmetic that are helpful for our purposes, nonetheless, for more on logarithmic arithmetic we urge the reader to consult the seminal paper of Jaulent \cite{Jaulent94}. Then we show some results on logarithmic ramification in \mbox{$\Zl$-extensions}. 

In Section 3 we construct a $\LL$-module, which is the logarithmic analogue to the $\LL$-module $X$ in \cite[\S 1]{Serre95}, to which we apply the codescent techniques. We show that this logarithmic $\LL$-module is noetherian and assuming the Gross-Kuz'min conjecture also torsion. 

In Section 4 we prove Theorem \ref{thm:IwasClog} in the non-cyclotomic case as well as a set of results involving the $\mut$ and $\lat$ logarithmic invariants for some class of $\Zl$-extensions.

Finally, in section 5 we provide numerical examples, in both the cyclotomic and non-cyclotomic case, of logarithmic class groups in the first layers of $\Zl$-extensions. Additionally we compare this information with the results in \cite{Hubbard&Washington17} to explicitely compute the $\mut$, $\lat$ and $\widetilde{\nu}$ logarithmic invariants. 

{\em Acknowledgements:} I am very grateful to professor Jean-Fran\c{c}ois Jaulent for having introduced me to the field of logarithmic arithmetic and his exemplary guidance. The author sincerely thanks the anonymous referee for the remarks and careful reading of the text, aiming to do the text more accessible.

\section{Logarithmic arithmetic}  
From now on we fix a number field $K$ and a prime number $\ell\neq 2$\footnote{The following results and statements have analogous formulations when $\ell=2$. In order to ease notations we will skip the treatment of such a prime.}. Let $\Pl_{K}^{0}$ be the set of finite places of $K$. \\

For any prime $p$ let $\widehat{\Qp^{c}}$ be the cyclotomic $\widehat{\Z}$-extension of $\Qp$. For a place $\p$ above $p$ we define the logarithmic ramification index $\widetilde{e}_{\p}$ and the logarithmic inertia degree $\widetilde{f}_{\p}$ as
\begin{equation*}
\widetilde{e}_{\p}=[\Kp:\Kp\cap\widehat{\Qp^{c}}] \	\	\	\text{ and }\	\	\	\widetilde{f}_{\p}=[\Kp\cap\widehat{\Qp^{c}}:\Qp]
\end{equation*}
respectively. Similarly, for a local extension $\Lp/\Kp$ we define the relative logarithmic ramification index $\widetilde{e}_{\pp|\p}$ and relative logarithmic inertia degree $\widetilde{f}_{\pp|\p}$ as
\begin{equation*}
\widetilde{e}_{\pp|\p}=[\Lp : \Lp\cap\widehat{\Q}_{p}^{\textrm{c}}\Kp]\	\	\	\text{ and }\	\	\	\	\widetilde{f}_{\pp|\p}=[\Lp\cap \widehat{\Q}_{p}^{\textrm{c}}\Kp:\Kp].
\end{equation*}

\begin{proposition}\label{prop:ramprop} Let $\Lp/\Kp$ be a local extension. The logarithmic ramification index and the logarithmic inertia degree satisfy the following properties:
\begin{enumerate}
\item Multiplicative relation: $n_{\p}=[\Kp:\Qp]=e_{\p}f_{\p}=\widetilde{e_{\p}}\widetilde{f_{\p}}$, where $e_{\p}$ and $f_{\p}$ are the usual ramification index and inertia degree.
\item For all $q \neq p$: $v_{q}(e_{\p})=v_{q}(\widetilde{e}_{\p})$.
\item Multiplicative in towers: $\widetilde{e}_{\mathfrak{P}}=\widetilde{e}_{\pp|\p}\widetilde{e}_{\p}$ and	$\widetilde{f}_{\mathfrak{P}}=\widetilde{f}_{\pp|\p}\widetilde{f}_{\p}$.
\end{enumerate}
\end{proposition}
\textit{Proof:} This is \cite[Thm. 1.4]{Jaulent94} and subsequent remarks. \hfill \qed

We endow $K$ with a family of $\ell$-adic logarithmic valuations $(\lval{\p})_{\p\in\Pl_{K}^{0}}$ taking values in $\Zl$ and which are defined as follows
$$\lval{\p}(x)=\left\lbrace \begin{array}{cl}
v_{\p}(x)& \textrm{ if }\p\nmid \ell, \\ & \\ \frac{\Log_{\ell}(\Norm_{\Kp/\Ql}(x))}{\deg \p }& \textrm{ if }\p\;|\;\ell,
\end{array}\right. \textrm{ for } x\in K_{\p}^{\times}.$$
Here $v_{\p}$ denotes the classical valuation, $\Norm_{\Kp/\Ql}$ is the local norm, $\Log_{\ell}$ is Iwasawa's logarithm and $\deg\p$ is a normalizing term defined\footnote{The definition of $\deg(\ell)$ given here differs to the one given in \cite[Def. 1.1]{Jaulent94}. Both definitions differ by a unit in $\Zl$ and this does not affect subsequent results. With this definition it is clearer that there is an element with logarithmic valuation 1 as shown in Remark \ref{rem:hp}.} as follows
$$\deg \p = \widetilde{f}_{\p} \deg p \	\	\	\textrm{ where }\	\	
\deg p=\left\lbrace
\begin{array}{cl} 
\Log_{\ell}(p) & \textrm{ if }p\neq\ell,  \\ \Log_{\ell}(1+\ell) & \textrm{ if }p=\ell.
   \end{array}\right.$$

\begin{remark}\label{rem:hp} It is clear that $\lval{\p}(x)$ takes values in $\Zl$ for $\p\nmid \ell$. This is in fact the case also for places $\p$ over $\ell$. Consider the function 
\[h_{\p}(x)=\frac{\Log_{\ell}(\Norm_{\Kp/\Ql}(x))}{n_{\p}\deg \p }.\]
If $\Kp=\Ql$ we have $\im(h_{\ell})=\Zl$, since $\Ql^{\times}\isom\ell^{\Z}\cdot \mu_{\ell-1}\cdot U_{\ell}^{(1)}$ and $\deg(\ell)=\ell\cdot\eta$ where $\eta\in\Zl^{\times}$. In particular $h_{\ell}(1+\ell)=1$. For general $\Kp$ we have that $\im(h_{\p})$ is the $\Zl$-module $[\Kp:\Kp\cap \Ql^{c}]^{-1}\Zl$ \cite[Prop. 1.2]{Jaulent94}, moreover $\lval{\p}=\widetilde{e_{\p}}h_{\p}$ \cite[Thm. 1.4]{Jaulent94}. In particular this makes the logarithmic divisor map below well defined.
\end{remark}

We define the logarithmic divisors $\mathcal{D}\ell_{K}$ as the free $\Zl$-module $\mathcal{D}\ell_{K}=\bigoplus_{\p\in\Pl_{K}^{0}} \Zl\p$ and
we extend additively the mapping $\deg$ to all of $\mathcal{D}\ell_{K}$, that is, for a logarithmic divisor $\mathfrak{d}$:
$$\deg(\mathfrak{d})=\deg\left( \sum_{\p\in\Pl_{K}^{0}}a_{\p}\p \right)= \sum_{\p\in\Pl_{K}^{0}}a_{\p}\deg \p, \	\	\	\textrm{ with }a_{\p}\in \Zl.$$

Let us denote $\Rc_{K}=\Zl\otimes_{\Z}K^{\times}$. The logarithmic divisor map 
$$\widetilde{div}:\Rc_{K}\longrightarrow\mathcal{D}\ell_{K},\	\	\	\	\	 x \mapsto \sum_{\p\in\Pl_{K}^{0}} \widetilde{v}_{\p}(x)\p,$$
associates to each  element in the tensor product a logarithmic divisor by means of the logarithmic valuation. It turns out that the image of an element in $\Rc_{K}$ has degree 0 \cite[Prop.\& Def. 2.2]{Jaulent94}. The image $\mathcal{P}\ell_{K}$ of all these elements is the  subgroup of principal logarithmic divisors. Therefore it is natural to consider the following exact sequence:
$$0\longrightarrow \widetilde{\mathcal{E}}_{K} \longrightarrow \Rc_{K} \stackrel{\widetilde{div}}{\longrightarrow} \mathcal{D}\ell_{K} \stackrel{\deg}{\longrightarrow} \Clog{K}^{*}  \longrightarrow 0,$$
we call $\widetilde{\mathcal{E}}_{K}$ the group of logarithmic units (defined in the obvious way) and $\Clog{K}^{*}\isom\mathcal{D}\ell_{K}/\mathcal{P}\ell_{K}$ the group of logarithmic classes of arbitrary degree. 

\begin{definition} The logarithmic class group $\Clog{K}\subset\Clog{K}^{*}$ consists of the classes of degree 0.
\end{definition}

The group $\Clog{K}^{*}$ of logarithmic classes of arbitrary degree is isomorphic to the Galois group $\Gal(K^{\lc}/K)$ of the maximal abelian $\ell$-extension of $K$ which splits completely over the cyclotomic $\Zl$-extension $K^{c}$ of $K$. The logarithmic class group is isomorphic to the relative Galois group $\Gal(K^{\lc}/K^{\text{c}})$ \cite[Thm. \& Conj. 2.3]{Jaulent94}. The Gross-Kuz'min conjecture asserts that the Galois group $\Gal(K^{\lc}/K)$ is a $\Zl$-module of rank 1, equivalently, asserts the finiteness of the logarithmic class group. 

\subsubsection*{Logarithmic ramification}

Let $\p$ be a finite place of $K$. Let $L/K$ be a finite extension, we say that $\p$ is logarithmically unramified if $\tilde{e}_{\pp|\p}=1$ for every place $\mathfrak{P}$ above $\p$. Hence if $L/K$ is infinite, we say that $\p$ is logarithmically unramified if and only if it is logarithmically unramified for every finite extension $K'/K$ contained in $L$. In the Galois case it might be more convenient to work with inertia groups to study ramification, so we give the corresponding description in the logarithmic setting in the very special case when $L$ is an abelian $\ell$-extension of $K$.

The so called $\ell$-adic class field theory \cite[I.1.2]{Jaulent86} establishes a correspondence between the decomposition group $D_{\p}$ of the maximal abelian $\ell$-extension $K^{\ab}$ of  $K$ with the pro-$\ell$-completion $\Rc_{\p}$ of the groups $K_{\p}^{\times}$, i.e. $D_{\p}\isom \Rc_{\p}= \varprojlim \Kp^{\times}/ \Kp^{\times \ell^{k}}$. We can extend the logarithmic valuations $\tilde{v}_{\p}$ to the $\Rc_{\p}$; in the sense that we have an analogous decomposition\footnote{Consider the short exact sequence of $\Zl$-modules $1\rightarrow \widetilde{\Uc}_{\p} \rightarrow \Rc_{\p} \stackrel{\lval{\p}}{\rightarrow} \Zl \rightarrow 0$, where the logarithmic valuation $\lval{\p}$ is a $\Zl$-epimorphism by \cite[Thm 1.4]{Jaulent94}. Then $\Rc_{\p}$ is the direct sum of $\widetilde{\Uc}_{\p}$ and a free $\Zl$-module of rank 1.} $\Rc_{\p}=\tilde{\pi}_{\p}^{\Zl}\widetilde{\Uc}_{\p}$ to the classical decomposition $\Rc_{\p}=\pi_{\p}^{\Zl}\Uc_{\p}$, where $\tilde{\pi}_{\p}$ is a logarithmic uniformizer and $\widetilde{\Uc}_{\p}$ is the subgroup of logarithmic units, i.e. the kernel under the corresponding logarithmic valuation. The image $\tilde{I}_{\p}\subset \Gal(K^{ab}/K)$ of $\widetilde{\Uc}_{\p}$ is said to be the subgroup of logarithmic inertia. Hence, we define the group of logarithmic inertia  $\tilde{I}_{\p}(L/K)$ of the abelian $\ell$-extension $L/K$ as the image of $\widetilde{\Uc}_{\p}$ in the Galois group $\Gal(L/K)$. We say $\p$ is logarithmically unramified in $L/K$ if $\tilde{I}_{\p}\subseteq \Gal(K^{\ab}/L)$. Notice that all the extensions we are considering are abelian, so we do not specify the place $\mathfrak{P}$ above $\p$ since the logarithmic inertia subgroups are isomorphic.

In general, $\ell$-adic class field theory (e.g. \cite[Thm. 1.1.13.]{Jaulent86}, \cite[\S 2]{Jaulent94}) establishes a correspondence between the Galois group of the maximal abelian $\ell$-extension of $K$ and the following quotient:
$$\Gal(K^{\ab}/K)\isom \Jc_{K}/\Rc_{K},$$
where $\Jc_{K}$ is the restricted product $\prod\Rc_{\p}$ of the $\ell$-adic completions $\Rc_{\p}$ defined above.

\begin{remark} Notice that in the finite case this characterisation of logarithmic \mbox{ramification} for abelian $\ell$-extensions in terms of inertia subgroups, coincides with the definition given at the beginning of the section. 

In fact, as we have seen, if $L/K$ is a finite abelian $\ell$-extension and $\mathfrak{P}$ is a place of $L$ above $\p\in\Pl_{K}^{0}$, the logarithmic inertia subgroup $\tilde{I}_{\p}(L/K)$ of the extension $L/K$ is the image of $\widetilde{\Uc}_{\p}$ inside the Galois group $\Gal(L/K)\isom\Jc_{K}/\Norm_{L/K}(\Jc_{L})\Rc_{K}$, hence:
$$\tilde{I}_{\mathfrak{p}}(L/K)\isom \widetilde{\Uc}_{\p}/\widetilde{\Uc}_{\p}\cap\Norm_{L/K}(\Jc_{L})\Rc_{K};$$
on the other hand, the logarithmic ramification index is defined as the degree $\widetilde{e}_{\pp|\p}=[\Lp:\Lp\cap\Kp^{c}]$, which is equal to the order of the Galois group $$\Gal(\Lp/\Lp\cap\Kp^{c})=\widetilde{\Uc}_{\p}\Norm_{\Lp/\Kp}(\Rc_{\mathfrak{P}})/\Norm_{\Lp/\Kp}(\Rc_{\mathfrak{P}}).$$

Since $L/K$ is abelian (\cite[Ch. 1 \S 2.a]{Jaulent86}) we have
$$\Norm_{\Lp/\Kp}(\Rc_{\mathfrak{P}})=\Rc_{\p}\cap\Norm_{L/K}(\Jc_{L})\Rc_{K}$$  
and so $\tilde{I}_{\mathfrak{p}}(L/K)\isom\Gal(\Lp/\Lp\cap\Kp^{c})$.
\end{remark}

Let $K_{\infty}$ be a $\Zl$-extension of $K$. We denote $K_{n}$ the $n$-th layer of $K_{\infty}$, i.e. $K_{n}$ is the unique subfield of $K_{\infty}$ with $\Gal(K_{n}/K)\isom \Z/\ell^{n}\Z$. 

As a consequence of the classic case, logarithmic ramification behaves well at the tame places in $\Zl$-extensions.

\begin{proposition} $\Zl$-extensions are logarithmically unramified at the places $\p\nmid\ell$ of $K$.\end{proposition}
\textit{Proof:} Let $K_{\infty}/K$ be a $\Zl$-extension. Let $\p$ be a place of $K$ over $p$ with $p\neq\ell$.  Let $\pp$ be a place in $K_{n}$ lying over $\p$. 

For every $n$ we have that $K_{n}$ is unramified at $\p$ \cite[Prop. 13.2]{Washington97}, i.e. $e_{\pp|\p}=1$. We have that $v_{q}(e_{\pp|\p})=v_{q}(\widetilde{e}_{\pp|\p})=0$ for every $q\neq p$, for the completions $K_{n,\pp}$ and $K_{\p}$ are both extensions of $\Qp$, hence the $q$-parts of the relative indices coincide due to (\textit{2}) and (\textit{3}) in Proposition \ref{prop:ramprop}. Moreover $\widetilde{e}_{\pp|\p}$ has no $p$-part since $p\nmid |\Gal(K_{n}/K)|$. We conclude that $\widetilde{e}_{\pp|\p}=1$ as desired.
\hfill \qed\\

In terms of ramification, the locally cyclotomic extension $K^{\lc}$ of $K$ corresponds to the maximal abelian extension of $K$ which is logarithmically unramified. Unlike the classical case, this extension $K^{\lc}$ is infinite since it contains the cyclotomic $\Zl$-extension $K^{\text{c}}$ of $K$. Since $K^{\lc}$ is a subextension of the maximal abelian $\ell$-ramified $\ell$-extension of $K$, the group $\Clog{K}^{*}\isom\Gal(K^{\lc}/K)$ is a noetherian $\Zl$-module, and assuming the Gross-Kuz'min Conjecture it has rank 1.  This yields the following result.

\begin{lemma}\label{lemma:RamLogauneplace} Let $K$ be a number field verifying the Gross-Kuz'min conjecture for the prime $\ell$, let $K_{\infty}/K$ be a non-cyclotomic $\Zl$-extension. Then at least one place $\p$ above $\ell$ is logarithmically ramified in $K_{\infty}/K$. 
\end{lemma}
\textit{Proof:} Let us suppose that $K_{\infty}$ is logarithmically unramified, hence it is contained in $K^{\lc}$, but then $\Gal(K^{\lc}/K)$ would have $\Zl$-rank two which contradicts the Gross-Kuz'min conjecture. \hfill \qed

\begin{remark} The cyclotomic $\Zl$-extension $K^{\text{c}}$ of $K$ is logarithmically unramified.
\end{remark}

Under the same assumptions and as a consequence of the preceding lemma we have the following corollary.

\begin{corollary}\label{cor:Existedlogtotram} There is an integer $d$ such that $K_{\infty}/K_{d}$ is logarithmically totally ramified for a non-empty subset of places above $\ell$ and logarithmically unramified (possibly inert) for the rest of the places.
\end{corollary}
\textit{Proof:} Let $\p_{1},\ldots,\p_{s}$ be the places in $K$ above $\ell$ which are logarithmically ramified. The logarithmic inertia subgroups $\tilde{I}_{i}$ are closed and non zero in $\Gamma=\Gal(K_{\infty}/K)$ ($i=1,\ldots,s$), therefore we have
$$\tilde{I}_{i}\isom \Gamma^{\ell^{d_{i}}},$$
for some $d_{i}\in\N$ ($i=1,\ldots,s$), taking $d=\max d_{i}$, it is clear that the extension $K_{\infty}/K_{d}$ is logarithmically totally ramified at every logarithmically ramified place. \hfill \qed

\begin{lemma} For a $\Zl$-extension $K_{\infty}$ of $K$ the following conditions are equivalent:
\begin{itemize}
\item[(i)] There exists a place $\p$ of $K$ such that $\widetilde{e}_{\pp|\p}>1$ for some place $\pp|\p$ of the subextension $K_{1}$ of $K$.
\item[(ii)] $K_{1}/K$ is logarithmically ramified.
\item[(iii)] $K_{\infty}$ is globally disjoint from the locally cyclotomic extension $K^{\lc}$ of $K$. 
\end{itemize}
\end{lemma}
\textit{Proof:} \textit{(i)}$\Leftrightarrow$ \textit{(ii)} is clear. \textit{(ii)} $\Leftrightarrow$ \textit{(iii)}: If $K_{\infty}\cap K^{\lc}=K_{n}$ for some $n>0$, then $K_{n}/K$ is logarithmically unramified. In particular $K_{1}/K$ is logarithmically unramified. 
\hfill \qed

\begin{definition} A $\Zl$-extension $K_{\infty}$ of $K$ is called locally disjoint from the cyclotomic $\Zl$-extension $K^{c}/K$ if it satisfies the above conditions.
\end{definition}

Let $K_{\infty}$ be a non-cyclotomic $\Zl$-extension of $K$ and $\p_{1},\ldots,\p_{s}$ be the places in $K$ above $\ell$ which ramify logarithmically. Define $n_{0}=\min \{d_{i}\;|\; \tilde{I}_{i}\isom \Gamma^{\ell^{d_{i}}}\}$. Then the $\Zl$-extension $K_{\infty}/K_{n_{0}}$ of $K_{n_{0}}$ is locally disjoint from the cyclotomic $\Zl$-extension $K_{n_{0}}^{\textrm{c}}/K_{n_{0}}$. 

\section{The logarithmic class group along non-cyclotomic $\Zl$-extensions}

Let $K$ be a number field admitting a non-cyclotomic $\Zl$-extension $K_{\infty}$ with Galois group $\Gamma=\Gal(K_{\infty}/K)$, and assume the Gross-Kuz'min conjecture is verified for $K$. We denote $K_{n}$ the unique subextension of $K_{\infty}$ such that $\Gal(K_{n}/K)\isom \Z/\ell^{n}\Z$. 

Let $\p_{1},\ldots,\p_{s}$ be the places of $K$ above $\ell$ which ramify logarithmically. There exist two integers $n_{0}$ and $d$ such that 
\begin{enumerate}
\item $0 \leq n_{0} \leq d$,
\item for every $n\geq n_{0}$ the $\Zl$-extension $K_{\infty}/K_{n}$ is locally disjoint from the cyclotomic  $\Zl$-extension $K^{\text{c}}_{n}=K_{n}K^{\text{c}}$ of $K_{n}$,
\item $K_{\infty}/K_{d}$ is completely logarithmically-ramified for some non-empty subset of wild places in $K_{d}$.
\end{enumerate}

Let $\widetilde{\mathcal{C}}^{*}$ be the inverse limit $\varprojlim \Clog{K_{n}}^{*}$ of the logarithmic classes of arbitrary degree of all the $K_{n}$. The group $\widetilde{\mathcal{C}}^{*}$ is isomorphic to the Galois group $\Gal(K_{\infty}^{\lc}/K_{\infty})$ of the maximal abelian logarithmically unramified extension of $K_{\infty}$. Hence $K_{n}^{\lc}$ is the maximal sub-extension of $K_{\infty}^{\lc}$ which is abelian over $K_{n}$ and splits completely over $K_{n}^{\textrm{c}}$, equivalently $K_{n}^{\lc}$ is the maximal sub-extension of $K_{\infty}^{\lc}$ which is abelian and logarithmically unramified over $K_{n}$. Hence for $n\geq d$ we have the following isomorphism of groups
\begin{equation*}
\Clog{K_{n}}^{*}\isom\Gal(K_{n}^{\lc}/K_{n})\isom\Gal(K_{\infty}K_{n}^{\lc}/K_{\infty}).
\end{equation*}

The profinite group $\Gamma$ acts on $\widetilde{\mathcal{C}}^{*}$ by innner automorphisms, since $\Gamma$ is procyclic the action is uniquely defined by a topological generator $\gamma\in\Gamma$. The group $G=\Gal(K_{\infty}^{\lc}/K)$ is then the semidirect product $\Gamma\ltimes\widetilde{\mathcal{C}}^{*}$. We can extend the action of $\Gamma$ to the Iwasawa algebra $\LL=\Zl[[\Gamma]]$ of $\Gamma$, making of $\widetilde{\mathcal{C}}^{*}$ a $\LL$-module. 

For every $n\geq 0$, we set
\begin{equation}\label{eq:cyclotomicfac}
\w_{n}=\gamma^{\ell^{n}}-1.
\end{equation}
We call cyclotomic factors the term $\w_{0}$ and the quotients $\w_{n}/\w_{n-1}$ for $n\geq 1$.

\begin{proposition}\label{prop:CetNoetandTors}
$\widetilde{\mathcal{C}}^{*}$ is a noetherian $\LL$-module. Moreover, if the Gross-Kuz'min conjecture holds for $\ell$ and every $K_{n}$ in $K_{\infty}$, it is torsion.
\end{proposition}
\textit{Proof:} The locally cyclotomic extension $K_{\infty}^{\lc}$ of $K_{\infty}$ is contained in the maximal abelian $\ell$-extension $K_{\infty}^{\lr}$ of $K_{\infty}$ which is $\ell$-ramified. The $\LL$-module $\widetilde{\mathcal{C}}^{*}$ can be written as a quotient of the noetherian $\LL$-module $X=\Gal(K_{\infty}^{\lr}/K_{\infty})$ \cite[Thm. 4]{Iwasawa73}. This yields the noetherian property. 

For $n\geq0$ let us denote $L_{\infty/n}$ the maximal abelian extension of $K_{n}$ which is logarithmically unramified over $K_{\infty}$, i.e. the maximal abelian extension of $K_{n}$ contained in $K^{\lc}_{\infty}$. Then
$$K_{\infty}\subseteq L_{\infty/0} \subseteq  L_{\infty/1} \subseteq \cdots \subseteq L_{\infty/n} \subseteq \cdots \subseteq K_{\infty}^{\lc},$$
\noindent and we have $K_{\infty}^{\lc}=\bigcup_{n\geq 0} L_{\infty/n}$. Then we have the following group isomorphisms
\[\w_{n}\widetilde{\Cc}^{*}\isom \Gal(K_{\infty}^{\lc}/L_{\infty/n})\	\	\	\	\	\widetilde{\Cc}^{*}/\w_{n}\widetilde{\Cc}^{*}\isom\Gal(L_{\infty/n}/K_{\infty}),\	\	\	n\geq 0.\]
\\
Let $d$ be as in Corollary \ref{cor:Existedlogtotram}, for $n \geq d$ let $\p_{1},\ldots,\p_{s}$ be the primes of $K_{n}$ which are logarithmically ramified in $K_{\infty}$; the number of these primes is the same for all $n\geq d$. For each $i=1,\ldots,s$ let $\tilde{I}_{i}$ be the logarithmic inertia subgroup of $\p_{i}$ in the extension $L_{\infty/n}/K_{n}$. Since $\p_{i}$ is totally logarithmically ramified in $K_{\infty}$ and the extension $L_{\infty/n}/K_{\infty}$ is logarithmically unramified, the groups $\tilde{I}_{i}$ are pro-cyclic isomorphic to the free $\Zl$-module $\Gamma^{\ell^{n}}\isom \Zl$. So we restrict $\tilde{I}_{i}$ to $\Gamma^{\ell^{n}}$. On the other hand, since none of the places apart of $\p_{1},\ldots,\p_{s}$ are logarithmically ramified in $L_{\infty/n}$, it follows that the product
$\tilde{I}_{1}\tilde{I}_{2}\cdots\tilde{I}_{s}$
fixes the maximal logarithmically unramified subextension of $L_{\infty/n}$ over $K_{n}$, which is nothing else than $K_{n}^{\lc}$, i.e. $\Gal(L_{\infty/n}/K_{n}^{\lc})=\tilde{I}_{1}\tilde{I}_{2}\cdots\tilde{I}_{s}$. Assuming the Gross-Kuz'min conjecture, $\Clog{K_{n}}^{*}\isom\Gal(K_{n}^{\lc}/K_{n})$ is a noetherian $\Zl$-module of rank 1. Then
\begin{eqnarray*}
\rank_{\Zl}\Gal(L_{\infty/n}/K_{n})&=&\rank_{\Zl}\Gal(K_{n}^{\lc}/K_{n})+\rank_{\Zl} \Gal(L_{\infty/n}/K_{n}^{\lc}) \\
&\leq & 1+s
\end{eqnarray*}
from which follows that
$$\rank_{\Zl}\Gal(L_{\infty/n}/K_{\infty})=\rank_{\Zl}\Gal(L_{\infty/n}/K)-1\leq s,$$
so
$\rank_{\Zl} \widetilde{\Cc}^{*}/\w_{n}\widetilde{\Cc}^{*}\leq s$,
for all $n\geq d$, therefore for all $n\geq 0$. Hence $\widetilde{\Cc}^{*}$ is $\LL$-torsion by \cite[Lemma 3 (ii)]{Iwasawa73}.
\hfill \qed

\begin{remark} We have just mimicked the classical proof of Iwasawa as in \cite[Thm. 5]{Iwasawa73}. In our case due to the structure of $\Clog{K_{n}}^{*}$ we have that 
$$\rank_{\Zl} \widetilde{\Cc}^{*}/\w_{n}\widetilde{\Cc}^{*} \leq s, \	\	\	\text{ for all }n\geq 0.$$
\end{remark}

\begin{proposition}\label{prop:CetisoCplusZl} Assume the Gross-Kuz'min conjecture for all $K_{n}$ in $K_{\infty}$. Let $\widetilde{\Cc}=\varprojlim\Clog{K_{n}}$ the inverse limit of the logarithmic class groups $\Clog{K_{n}}$. There is an isomorphism of $\Zl$-modules
\begin{equation}\label{isom:CetisoCplusZl}
\widetilde{\Cc}^{*}\isom \widetilde{\Cc}\oplus \Zl.
\end{equation}
\end{proposition}
\textit{Proof:} Consider the following exact sequence
\begin{equation*}
0\rightarrow \Clog{K_{n}} \rightarrow \Clog{K_{n}}^{*} \rightarrow \Zl \rightarrow 0.
\end{equation*}
We take the inverse limit $\varprojlim$ for the arithmetic norm maps in each of the terms. It is a well known fact that profinite groups satisfy the Mittag-Leffler property (e.g. \cite[II.7]{Neukirch&Schmidt&Wingberg08}). So we obtain exactness when taking inverse limits
\begin{equation*}
0\rightarrow \widetilde{\Cc} \rightarrow \widetilde{\Cc}^{*} \rightarrow \Zl \rightarrow 0.
\end{equation*}
 \hfill \qed

\subsubsection*{Assumption on the logarithmic ramification} 

Before approaching our problem in great generality, let us assume first that the places of $K$ which ramify logarithmically in $K_{\infty}$ are totally logarithmically ramified. We can always do so by rewriting the $\Zl$-extension in order to have $K=K_{d}$, with $d$ as in Corollary \ref{cor:Existedlogtotram}.

For each place $\p_{i}$ ramifying logarithmically, let $G=\Gal(K^{\lc}_{\infty}/K)\isom \widetilde{\Cc}^{*}\rtimes\Gamma$ and let $\tilde{I}_{i}\subset G$ its logarithmic inertia subgroup. Since $K_{\infty}^{\lc}/K_{\infty}$ is logarithmically unramified, we have that each logarithmic inertia subgroup has trivial intersection with the group
$\Gal(K_{\infty}^{\lc}/K_{\infty})$: $\tilde{I}_{i}\cap \Gal(K_{\infty}^{\lc}/K_{\infty})=1$. Since $K_{\infty}/K$ is totally logarithmically ramified over $\p_{i}$, the groups of logarithmic inertia $\tilde{I}_{i}$ restrict to $\Gamma$.

Since we are mainly interested in the study of the logarithmic class groups (of zero degree) $\Clog{K_{n}}$, we consider the group $H=\Gal(K_{\infty}^{\lc}/K^{c})$. We take advantage of the fact that the cyclotomic extension is logarithmically unramified, so the logarithmic ramification groups $\tilde{I}_{i}\subset G$ remain unchanged in $H$, i.e. they are procyclic and restrict to $\Gamma$. The group $H$ is the semi-direct product of $\Gamma$ and the normal subgroup $\Gal(K^{\lc}_{\infty}/K_{\infty})$, and can be written as
$$H=\tilde{I}_{i}\Gal(K^{\lc}_{\infty}/K_{\infty})=\Gal(K^{\lc}_{\infty}/K_{\infty})\tilde{I}_{i}\	\	\text{ for }\	i=1,\ldots,s.$$
Fixing topological generators $\tilde{\sigma}_{i}$ of $\tilde{I}_{i}$, there exists $\tilde{a}_{1},\ldots,\tilde{a}_{s}\in \Gal(K^{\lc}_{\infty}/K_{\infty})$ such that 
$$\tilde{\sigma}_{i}=\tilde{a}_{i}\tilde{\sigma}_{1}.$$ 
In particular we can take $\tilde{\sigma}_{1}=\gamma$ for a topological generator $\gamma$ of $\Gamma$.

Let $H'$ be the closure of the commutator subgroup of $H$. We know that $H'\isom\w_{0}\widetilde{\mathcal{C}}\isom T\widetilde{\mathcal{C}}$, the last is the isomorphism $\Zl[[\Gamma]\isom \Zl[[T]]$ given by $\gamma \mapsto T+1$ \cite[Lemma 13.14]{Washington97}. 

\begin{remark} After what we have been discussing we have the next result whose proof goes the same way as in the classical case of \cite[Lemma 13.15]{Washington97}. Nonetheless we reproduce the proof in our special setting. 
\end{remark}

\begin{proposition}\label{prop:GCLisoCtildQuot}
Let $Y_{0}$ be the $\Zl$-submodule of $\widetilde{\mathcal{C}}$ generated by $\{\tilde{a}_{i}\,|\,2\leq i \leq s\}$ and $\w_{0}\widetilde{\mathcal{C}}=T\widetilde{\mathcal{C}}$. Let $Y_{n}=\frac{\w_{n}}{\w_{0}}Y_{0}$, where
$$\frac{\w_{n}}{\w_{0}}=1+\gamma+\gamma^{2}+\cdots +\gamma^{\ell^{n}-1}=\frac{(1+T)^{\ell^{n}}-1}{T}.$$
Then
$$\Clog{K_{n}}\isom \widetilde{\mathcal{C}}/Y_{n} \	\	\textrm{ for }n\geq 0.$$
\end{proposition}
\textit{Proof:} First, we show the case $\Clog{K}\isom \widetilde{\mathcal{C}}/Y_{0}$. The locally cyclotomic extension $K^{\lc}$ is the maximal abelian extension logarithmically unramified of $K^{\textrm{c}}$ contained in $K^{\lc}_{\infty}$. Its Galois group $\Gal(K^{\lc}/K^{\textrm{c}})$ is hence equal to $\Gal(K^{\lc}_{\infty}/K^{\textrm{c}})/Z_{0}$, where $Z_{0}$ is the smallest subgroup of $H=\Gal(K^{\lc}_{\infty}/K^{\text{c}})$ containing the commutator subgroup $H'$ as well as the logarithmic inertia subgroups $\tilde{I}_{i}$, $1\leq i \leq s$. The decomposition of $H$ in a semi-direct product of $\tilde{I}_{1}$ and $\widetilde{\mathcal{C}}$ yields
\begin{eqnarray*}
\Clog{K} \isom \Gal(K^{\lc}/K^{\textrm{c}})  & \isom & H/\Gal(K_{\infty}^{\lc}/K^{\lc})\\
& \isom & \widetilde{\mathcal{C}}\tilde{I}_{1}/\overline{\langle\tilde{I}_{1},\tilde{a}_{2}\ldots,\tilde{a}_{s},\w_{0}\widetilde{\mathcal{C}}\rangle}\\
& \isom & \widetilde{\mathcal{C}}/\overline{\langle\tilde{a}_{2},\ldots,\tilde{a}_{s},\w_{0}\widetilde{\mathcal{C}}\rangle} =  \widetilde{\mathcal{C}}/Y_{0}.
\end{eqnarray*}

Now let us assume $n\geq 1$. If we replace $K^{\textrm{c}}$ by $K_{n}^{\textrm{c}}$, the extension $K_{\infty}^{c}/K_{n}^{\textrm{c}}$ is also a $\Zl$-extension. The above result can be applied by replacing the topological generator $\gamma$ by $\gamma^{\ell^{\scriptscriptstyle \ssp{\rule{0ex}{1ex}}{n}}}$ and the $\tilde{\sigma_{i}}$ become $\tilde{\sigma}_{i}^{\ell^{\scriptscriptstyle \ssp{\rule{0ex}{1ex}}{n}}}$. Notice that for every integer $k$, we have:
\begin{eqnarray*}
\tilde{\sigma}_{i}^{k}&=&(\tilde{a}_{i}\tilde{\sigma}_{1})^{k}\\
&=& \tilde{a}_{i}
(\tilde{\sigma}_{1}\tilde{a}_{i}\tilde{\sigma}_{1}^{-1})
(\tilde{\sigma}_{1}^{2}\tilde{a}_{i}\tilde{\sigma}_{1}^{-2})\cdots
(\tilde{\sigma}_{1}^{k-1}\tilde{a}_{i}\tilde{\sigma}_{1}^{-(k-1)})
\tilde{\sigma}_{1}^{k}\\
&=& \tilde{a}_{i}^{1+\tilde{\sigma}_{1}+\ldots +\tilde{\sigma}_{1}^{k-1}}\tilde{\sigma}_{1}^{k}
\end{eqnarray*}
and particularly, for $k=\ell^{n}$, the isomorphism $\tilde{I}_{1}\isom\Gamma$ gives $\tilde{\sigma}_{i}^{\ell^{n}}=\tilde{a}_{i}^{\w_{n}/\w_{0}}(\tilde{\sigma}_{1})^{\ell^{n}}$. This implies that the groups of logarithmic inertia of $H_{n}=\Gal(K_{\infty}^{\lc}/K_{n}^{\textrm{c}})$ are generated by the $\frac{\w_{n}}{\w_{0}}\tilde{a}_{i}$. Finally, the commutator subgroup $H'_{n}$ of $H_{n}$ is given by $\omega_{n}\widetilde{\mathcal{C}}$.
From which the result follows.   \hfill \qed

\subsubsection*{General case}

In what follows we come back to the general situation exposed in the beginning of the section. Let $K_{\infty}$ be a non-cyclotomic $\Zl$-extension of $K$. Let $\p_{1},\ldots,\p_{s}$ be the places of $K$ which ramify logarithmically and let $d$ the smallest integer such that the $\Zl$-extension $K_{\infty}/K_{d}$ is totally logarithmically ramified.  

\begin{proposition}\label{prop:dsuchthatClogisomquoti}
Let $K_{\infty}$ be a non-cyclotomic $\Zl$-extension of $K$. There exists $d \geq 0$ such that
$$\Clog{K_{n}}\isom \widetilde{\Cc}/Y_{n}\	\	\	\textrm{ for }\	\	n \geq d.$$
\end{proposition} 
\textit{Proof:} We apply the Proposition \ref{prop:GCLisoCtildQuot} to the extension $K_{\infty}/K_{d}$, the group $\widetilde{\Cc}$ being the same as the original extension. The element $\gamma^{\ell^{d}}$ generates $\Gal(K_{\infty}/K_{d}^{\textrm{c}})$, and naturally $\omega_{d}\widetilde{\Cc}$ generates the commutator subgroup of $\Gal(K_{\infty}/K_{d})$. For $n\geq d$ we have
$$\frac{\w_{n}}{\w_{d}}=1+\gamma^{\ell^{d}}+\gamma^{2\ell^{d}}+\cdots +\gamma^{\ell^{n}-\ell^{d}},$$
it comes that the commutator subgroup of $\Gal(K_{\infty}/K_{n})$ is given by $\frac{\w_{n}}{\w_{d}}(\w_{d}\widetilde{\Cc})=\w_{n}\widetilde{\Cc};$ in general for the $Y_{n}$ we have
$Y_{n}=\frac{\w_{n}}{\w_{d}}Y_{d}$.

\hfill \qed
\\
The exact sequence
\begin{equation} \label{eq:pseudo}
0\rightarrow Y_{d} \rightarrow \widetilde{\Cc} \rightarrow \Clog{K_{d}} \rightarrow 0
\end{equation} 
sets a pseudo-isomorphism $Y_{d}\sim \widetilde{\Cc}$, that is a morphism $Y_{d} \rightarrow \widetilde{\Cc}$ of $\LL$-modules with finite kernel and cokernel. 

\section{Structural parameters}

Let us recall a fundamental result (e.g. \cite[Thm 13.12]{Washington97}).

\begin{theorem}[Structure theorem of noetherian $\LL$-modules]  \label{thm:structure}
Let $M$ be a noetherian $\LL$-module. There is a $\LL$-module morphism $M\rightarrow E$ to a unique elementary $\LL$-module $E$ with finite kernel and finite cokernel:
$$M \sim E=\LL^{\rho}
\oplus\left(\bigoplus\limits_{i=1}^{s}\LL/\ell^{m_{i}}\LL\right)
\oplus\left(\bigoplus\limits_{j=1}^{t}\LL/P_{j}\LL\right)$$
with $P_{j}$ distinguished polynomials, ordered by divisibility: $P_{1}|P_{2}|\ldots|P_{t}$. The integer $\rho$ is the rank of the $\LL$-module $M$, and the polynomial
$$\x=\prod_{i=1}^{s}\ell^{m_{i}}\prod_{j=1}^{t}P_{j}$$
is the characteristic polynomial of the torsion sub-module of $M$. 
\end{theorem}

\begin{definition}\label{def:invariants}
With the same notation of the above theorem let's set
$$\rho = \rank_{\LL}M,\	\	\	\mu=\sum_{i=1}^{s} m_{i} \	\	\	\textrm{ and }\	\	\	\lambda=\sum_{j=1}^{t}\deg(P_{j});$$
we say that $\rho$, $\mu$ and $\lambda$ are the structural invariants of the $\LL$-module $M$.
\end{definition}

We have that $\widetilde{\Cc}$ is a noetherian $\LL$-module (Prop. \ref{prop:CetNoetandTors}). Let us write $\tilde{\rho}$, $\tilde{\mu}$ and $\tilde{\lambda}$ for the invariants defined above for $\widetilde{\Cc}$. Therefore assuming the Gross-Kuzmin conjecture we have $\tilde{\rho}=0$ and we obtain the following theorem.

\begin{theorem}\label{thm:Iwalog} Let $K_{\infty}$ be a non-cyclotomic $\Zl$-extension of $K$. Let $\Clog{K_{n}}$ be the logarithmic class group of $K_{n}$. Assuming the Gross-Kuz'min conjecture for $\ell$ and every finite subextension $K_{n}$ of $K_{\infty}/K$ there exist integers $\tilde{\lambda}, \tilde{\mu} \geq 0$ and $\tilde{\nu}$, such that
$$|\Clog{K_{n}}|=\ell^{\tilde{\mu}\ell^{n}+\tilde{\lambda} n + \tilde{\nu}},$$
for $n$ big enough. 
\end{theorem}

\begin{remark} The above theorem together with Jaulent's results in the cyclotomic case yield Theorem 1. 
\end{remark}

In order to compare the classical and logarithmic invariants, we recall the following definition.

\begin{definition} We say that a sequence $(X_{n})_{n\in\N}$ of finite $\Zl$-modules is parametrized by $(\rho,\mu,\lambda)$ if the order $\ell^{x(n)}$ of $X_{n}$ is asymptotically given by the formula:
\begin{eqnarray*}
x(n) & \approx &  (n+1)\rho{\ell^{n}}+ \mu\ell^{n}+  \lambda n,
\end{eqnarray*}
where $\approx$ means that the difference of the right and left side terms is bounded. Moreover, we say that it is strictly parametrized if the difference is eventually constant. 
\end{definition} 

Let $K_{\infty}/K$ be a non-cyclotomic $\Zl$-extension of $K$. Let $\Cc=\varprojlim C\ell$ be the inverse limit of the $\ell$-class groups of the $K_{n}$. Then $\Cc$ is a torsion $\LL$-module corresponding to the Galois group $\Gal(K_{\infty}^{\nr}/K_{\infty})$ of the maximal abelian unramified $\ell$-extension $K_{\infty}^{\nr}$ of $K_{\infty}$. Similarly the inverse limit $\widetilde{\Cc}^{*}=\varprojlim \Clog{K_{n}}^{*}$ of the logarithmic classes of arbitrary degree of the $K_{n}$, corresponds to the Galois group $\Gal(K_{\infty}^{\lc}/K_{\infty})$ of the maximal abelian logarithmically unramified $\ell$-extension $K_{\infty}^{\lc}$ of $K_{\infty}$; after the propositions \ref{prop:CetNoetandTors} and \ref{prop:CetisoCplusZl}, assuming the Gross-Kuz'min conjecture $\widetilde{\Cc}^{*}$ is a noetherian torsion $\LL$-module, isomorphic as $\Zl$-module to a direct sum $\widetilde{\Cc}^{*}\isom\widetilde{\Cc}\oplus\Zl$. As a consequence of theorem \ref{thm:structure} we obtain the following proposition.

\begin{proposition}\label{charpolinva} Let $\widetilde{\Cc}^{*}$ and $\widetilde{\Cc}$ be as before. The characteristic polynomials associated to these $\LL$-modules satisfy
\begin{equation*}
\x_{{}_{\widetilde{\Cc}^{*}}}=\x_{{}_{\widetilde{\Cc}}} \cdot (\gamma-1),
\end{equation*}
and we have the following relations for their parameters
\begin{equation}\label{eq:ParLogArbLogCero}
\widetilde{\mu}^{*} = \widetilde{\mu}\	\	\	\textrm{ and }\	\	\	\	\widetilde{\lambda}^{*}=\widetilde{\lambda}+1,
\end{equation}
where $\widetilde{\mu}^{*}$ and  $\widetilde{\lambda}^{*}$ (respectively $\widetilde{\mu}$ and  $\widetilde{\lambda}$) are the parameters associated to the $\LL$-module $\widetilde{\Cc}^{*}$ (respectively $\widetilde{\Cc}$).
\end{proposition}

Let $\p_{1},\ldots,\p_{s}$ be the places of $K$ above $\ell$. We have the short exact sequence of $\ell$-groups
\begin{equation}\label{eq:SuiteLogEt}
0 \longrightarrow \Clog{K}^{* [\ell]} \longrightarrow \Clog{K}^{*} \stackrel{\Theta}{\longrightarrow} C\ell'_{K} \longrightarrow 0,
\end{equation}
where $\Clog{K}^{* [\ell]}$ is the subgroup of $\Clog{K}^{*}$ generated by the classes of the  $\p_{i}$; the group $C\ell'_{K}$ is the quotient of the $\ell$-part of the class group of $K$ and the $\ell$-part of the subgroup generated by the classes associated to the primes $\p_{i}$ above $\ell$
\begin{equation}\label{eq:SuiteLogCel}
0 \longrightarrow C\ell^{[\ell]}_{K} \longrightarrow C\ell_{K} \longrightarrow C\ell'_{K} \longrightarrow 0;
\end{equation}
and $\Theta$ is defined as 
\begin{equation*}\label{eq:ThetaLog}
\Theta : \sum_{\p}m_{\p}\p \mapsto \prod_{\p\nmid\ell} \p^{m_{\p}}.
\end{equation*}

Since the Mittag-Leffler property holds on profinite groups, taking the inverse limit in the short exact sequences \eqref{eq:SuiteLogEt} and \eqref{eq:SuiteLogCel} gives short exact sequences, and setting $\Cc'= \varprojlim C\ell'_{K_{n}}$ the corresponding isomorphisms 
\begin{equation}\label{isom:CcprimeIsom}
\Cc' \isom \widetilde{\Cc}^{*}/\widetilde{\Cc}^{*[\ell]} \isom \Cc/\Cc^{[\ell]}.
\end{equation}

Le us denote by $\mu$ and $\lambda$ (respectively $\mu'$, $\lambda'$) the parameters attached to the $\LL$-module $\Cc$ (respectively $\Cc'$). We have the following obvious relations
\begin{equation}\label{eq:ParMuLa} \begin{array}{ccc}
\widetilde{\mu}^{*}=\widetilde{\mu}&\widetilde{\mu}^{*} = \mu' + \widetilde{\mu}^{*[\ell]}&\widetilde{\lambda}^{*}=\lambda' + \widetilde{\lambda}^{*[\ell]} \\ 
\widetilde{\lambda}^{*}=\widetilde{\lambda}+1&\mu  =\mu'+\mu^{[\ell]}&\lambda =\lambda'+\lambda^{[\ell]} 
\end{array}\end{equation}
with $\mu^{[\ell]}$ and $\lambda^{[\ell]}$ (respectively $\widetilde{\mu}^{*[\ell]}$ and $\widetilde{\lambda}^{*[\ell]}$) as the associated parameters to $\Cc^{[\ell]}$ (respectively $\widetilde{\Cc}^{*[\ell]}$).

The following definition, if verified, allows us to have some control on the characteristic polynomials of the preceding modules.

\begin{definition}\label{def:DecompoFinie}
Let $K_{\infty}$ be a $\Zl$-extension of $K$. We say that a place $\p$ of $K$ splits finitely in the tower $K_{\infty}/K$, if its decomposition group in $\Gal(K_{\infty}/K)$ is open, equivalently, if there is an integer $n_{\p}$ such that the layer $K_{n_{\p}}$ corresponds to the decomposition field of $\p$ in $K_{\infty}/K$.
\end{definition}

\begin{remark}
The preceding definition holds for every wild place in the cyclotomic $\Zl$-extension of every number field.
\end{remark}

With the notation of \eqref{eq:ParMuLa} and \eqref{eq:cyclotomicfac} we state the following theorem: 

\begin{theorem}\label{Thm:relinva} Let $K_{\infty}/K$ be a non-cyclotomic $\Zl$-extension. Assume that the places of $K$ above $\ell$  split finitely in the extension $K_{\infty}/K$. Then the structural invariants $(\rho,\mu,\lambda)$ attached to the sequence $(C\ell_{n})_{n\in\N}$ of $\ell$-class groups and the invariants $(\widetilde{\rho},\widetilde{\mu},\widetilde{\lambda})$ attached to the sequence $(\Clog{n})_{n\in\N}$ of logarithmic class groups along the extension $K_{\infty}/K$ satisfy the following relations:
$$\rho=\widetilde{\rho}=0,\		\	\	\	\mu=\widetilde{\mu}\		\	\	\	\textrm{ but }\	\	\	\	\lambda \textrm{ might be different of } \tilde{\lambda}.$$
And the characteristic polynomials of the torsion $\LL$-modules $\Cc$ and $\widetilde{\Cc}$ just differ by cyclotomic factors.
\end{theorem}
\textit{Proof:} The first relation is trivial since $\Cc$ and $\widetilde{\Cc}^{*}$ are $\LL$-torsion modules. 

Let $D_{i}$ be the decomposition subgroup for the place $\p_{i}$ in $K_{\infty}/K$; since the $D_{i}$ are open they are isomorphic to $\Gamma^{\ell^{m_{i}}}$ for some $m_{i}$. Take $m=\max_{1\leq i\leq s}\{ m_{i}\}$, then $\Gamma^{\ell^{m}}$ acts trivially on $\Cc^{[\ell]}$ and $\widetilde{\Cc}^{*[\ell]}$, therefore $\w_{m}=\gamma^{\ell^{m}}-1$ annihilates them, i.e. their characteristic polynomial has no $\ell$-power factors which implies $\mu^{[\ell]}=\widetilde{\mu}^{*[\ell]}=0$.

We have the decomposition into cyclotomic factors
$$\w_{m}=\w_{0}\frac{\w_{1}}{\w_{0}}\frac{\w_{2}}{\w_{1}}\cdots \frac{\w_{m}}{\w_{m-1}}.$$
Since the minimal polynomial and the characteristic polynomials must have the same factors, the characteristic polynomials of $\Cc^{[\ell]}$ and $\widetilde{\Cc}^{*[\ell]}$ differ just by cyclotomic factors. 

Hence the $\LL$-modules $\Cc$ and $\widetilde{\Cc}^{*}$ have the same $\mu$ invariants by \eqref{eq:ParMuLa} and their $\lambda$ invariants might be different. From \eqref{eq:ParLogArbLogCero} we recover the relations for the logarithmic class group. \hfill \qed

\begin{corollary} With the same notation as above we have
\[\mu = \widetilde{\mu}=\widetilde{\mu}^{*}=\mu'.\]
\end{corollary}

\begin{remark}  Assuming the Gross-Kuz'min conjecture in Proposition \ref{prop:CetisoCplusZl} we deduced that $\widetilde{\lambda}^{*}=\widetilde{\lambda}+1$. Observing that the $\Zl$-submodule $\widetilde{\Cc}^{*[\ell]}$ of $\widetilde{\Cc}^{*}$ is infinite, for it contains the classes $\alpha\cdot\p$ for $\p$ a place over $\ell$ and $\alpha\in\Zl$ (see \S 2), we can apply an analogous argument to the $\LL$-modules $\widetilde{\Cc}^{[\ell]}$ and $\widetilde{\Cc}^{*[\ell]}$, from which we deduce that $\widetilde{\lambda}^{*[\ell]}=\widetilde{\lambda}^{[\ell]}+1$.
\end{remark}

The following table summarizes our results. 

\begin{center}
\begin{tabular}{|c|c|c|c|c|}
\hline
  \multicolumn{5}{|c|}{Invariants attached to parametrized sequences of finite $\Zl$-modules.} \\
  \hline
  \multirow{2}{*}{Modules}     & \multirow{2}{*}{Definition} & \multicolumn{3}{c|}{Parameters} \\
  \cline{3-5}
    &                             & $\rho$ &     $\lambda$        & $\mu$ \\ 
  \hline
   \multirow{2}{*}{$C\ell'_{n}$} & \multirow{2}{*}{$C\ell'_{n}$ is the $\ell$-group of $\ell$-classes of $K_{n}$} & \multirow{2}{*}{0} & \multirow{2}{*}{$\lambda'$} & \multirow{2}{*}{$\mu'$} 
   \\
   & & & & 
   \\
   \hline
   \multirow{2}{*}{$\Clog{n}$} & \multirow{2}{*}{$\Clog{n}$ is the logarithmic class group of $K_{n}$} & \multirow{2}{*}{0} & \multirow{2}{*}{$\tilde{\lambda}=\lambda'+\widetilde{\lambda}^{[\ell]}$} & \multirow{2}{*}{$\mu'$} \\ 
   & & & & \\
   \hline
   \multirow{2}{*}{$C\ell_{n}$} & \multirow{2}{*}{$C\ell_{n}$ is the $\ell$-class group of $K_{n}$} & \multirow{2}{*}{0} & \multirow{2}{*}{$\lambda=\lambda'+\lambda^{[\ell]}$} & \multirow{2}{*}{$\mu'$} \\ 
   & & & & \\
  \hline
\end{tabular}
\end{center}

We give an example for which we can apply our results.

\begin{example} Let $K$ be a quadratic imaginary number field and $\ell$ a prime number splitting in $K$, say $(\ell)=\mathfrak{l}\bar{\mathfrak{l}}$. The maximal abelian  $\mathfrak{l}$-ramified $\ell$-extension $M_{\mathfrak{l}}$ of $K$ contains a non-cyclotomic $\Zl$-extension of $K$, say $K_{\infty}$. For consider the $\ell$-Hilbert class field $H$ of $K$, then 
\begin{eqnarray*}
\Gal(M_{\mathfrak{l}}/H) &=&\prod_{\p}\Uc_{\p}\Rc_{K}/\prod_{\p\neq \mathfrak{l}}\Uc_{\p}\Rc_{K} \\
&=&\Uc_{\mathfrak{l}}/\Uc_{\mathfrak{l}}\cap(\prod_{\p\neq \mathfrak{l}}\Uc_{\p}\Rc_{K}) \\
&=&\Uc_{\mathfrak{l}}/s_{\mathfrak{l}}(\mu_{K}) \\
&=&\Zl,
\end{eqnarray*}
where $s_{\mathfrak{l}}(\mu_{K})$ is the image of the generalized units $\mu_{K}$ of $K$ under the localisation map $s_{\mathfrak{l}}$ \cite[Thm. \& Def. 1.1.22]{Jaulent86}. Therefore $\Gal(M_{\mathfrak{l}}/K)\isom \Zl$ $\times$ finite abelian $\ell$-group. 

 In particular, since $K/\Q$ is abelian, Lemma  \ref{lemma:RamLogauneplace} states that at least one place $\mathfrak{l}$ over $\ell$ ramifies logarithmically in $K_{\infty}$. Finally, the image of the decomposition subgroup $D_{\bar{\mathfrak{l}}}$ inside $\Gal(M_{\mathfrak{l}}/K)$ is isomorphic to $\Zl$ $\times$ finite abelian $\ell$-group, hence it has finite index in $\Gal(M_{\mathfrak{l}}/K)$. For the image of $D_{\bar{\mathfrak{l}}}$ is
$$\Rc_{\bar{\mathfrak{l}}}/\Rc_{\bar{\mathfrak{l}}}\cap(\prod_{\p\neq \mathfrak{l}}\Uc_{\p}\Rc_{K}).$$
Therefore both places $\mathfrak{l}$ and $\bar{\mathfrak{l}}$ split finitely.
\end{example}

\section{Numerical examples} \label{sec:numericalex}

\subsection{Quadratic fields}\label{sub:quadfiels}

We use the functions and algorithms implemented by Belabas and Jaulent in \cite{Belabas&Jaulent16} to compute the logarithmic class group of a number field. These algorithms are now available in PARI/GP \cite{PARI2}. We use PARI/GP to compute the extensions corresponding to the first layers of some $\Zl$-extensions, afterwards we apply the logarithmic class group algorithms to each of these number fields.

Let $K$ be a quadratic field of the form $\Q(\sqrt{d})$. For now we will restrict to the computation of the 3-group of logarithmic classes in the layers $K=K_{0}$, $K_{1}$ and $K_{2}$ of the cyclotomic $\Z_{3}$-extension of $K$.

In the next table we take $-100 < d < 100$. For sake of space, we just show all the cases for which the logarithmic class group is not trivial and for which $3$ splits in $\Q(\sqrt{d})$.

\begin{center}
\begin{tabular}{c|ccc|ccc|ccc|}
\cline{2-10}
  & \multicolumn{3}{c|}{\multirow{2}{*}{$\Q(\sqrt{d})$}} & \multicolumn{3}{c|}{\multirow{2}{*}{$\Q(\sqrt{d},\cos(2\pi/9))$}} & \multicolumn{3}{c|}{\multirow{2}{*}{$\Q(\sqrt{d},\cos(2\pi/27))$}}\\
   & & & & & & & & & \\
  \hline
   \multicolumn{1}{ |c|}{\multirow{2}{*}{$d$}}     & \multicolumn{1}{ c|}{\multirow{2}{*}{$\Clog{K}$}} & \multicolumn{1}{ c|}{\multirow{2}{*}{$\Clog{K}^{[\ell]}$}} & \multirow{2}{*}{$C\ell'_{K}$} & \multicolumn{1}{ c|}{\multirow{2}{*}{$\Clog{K_{1}}$}} & \multicolumn{1}{ c|}{\multirow{2}{*}{$\Clog{K_{1}}^{[\ell]}$}} & \multirow{2}{*}{$C\ell'_{K_{1}}$} & \multirow{2}{*}{$\Clog{K_{2}}$} & \multicolumn{1}{ |c|}{\multirow{2}{*}{$\Clog{K_{2}}^{[\ell]}$}} & \multirow{2}{*}{$C\ell'_{K_{2}}$}  \\
   
    \multicolumn{1}{ |c|}{} & \multicolumn{1}{ c|}{} & \multicolumn{1}{ c|}{} & & \multicolumn{1}{ c|}{}  & \multicolumn{1}{ c|}{} & \multicolumn{1}{ c|}{} & & \multicolumn{1}{ |c|}{} & \\
  \hline
  \multicolumn{1}{|c|}{-86}  & \multicolumn{1}{ c|}{[3]} & [3] & \multicolumn{1}{ |c|}{[]}  & \multicolumn{1}{ c|}{[9, 9]} & \multicolumn{1}{ c|}{[3]} & [9, 3] &[27, 27] & \multicolumn{1}{ |c|}{[3]} & [27, 9]\\
       \hline 
\multicolumn{1}{|c|}{-74}  & \multicolumn{1}{ c|}{[9]} & [9] & \multicolumn{1}{ |c|}{[]}  & \multicolumn{1}{ c|}{[27]} & \multicolumn{1}{ c|}{[9]} & [3] &[81] & \multicolumn{1}{ |c|}{[9]} & [9]\\
       \hline 
\multicolumn{1}{|c|}{-65}  & \multicolumn{1}{ c|}{[3]} & [3] & \multicolumn{1}{ |c|}{[]}  & \multicolumn{1}{ c|}{[9]} & \multicolumn{1}{ c|}{[3]} & [3] &[27] & \multicolumn{1}{ |c|}{[3]} & [9]\\
       \hline 
\multicolumn{1}{|c|}{-47}  & \multicolumn{1}{ c|}{[9]} & [9] & \multicolumn{1}{ |c|}{[]}  & \multicolumn{1}{ c|}{[27]} & \multicolumn{1}{ c|}{[9]} & [3] &[81] & \multicolumn{1}{ |c|}{[9]} & [9]\\
       \hline 
\multicolumn{1}{|c|}{-41}  & \multicolumn{1}{ c|}{[27]} & [27] & \multicolumn{1}{ |c|}{[]}  & \multicolumn{1}{ c|}{[81, 3]} & \multicolumn{1}{ c|}{[27]} & [3, 3] &[243, 9] & \multicolumn{1}{ |c|}{[27]} & [9, 9]\\
       \hline 
\multicolumn{1}{|c|}{-35}  & \multicolumn{1}{ c|}{[3]} & [3] & \multicolumn{1}{ |c|}{[]}  & \multicolumn{1}{ c|}{[9]} & \multicolumn{1}{ c|}{[3]} & [3] &[27] & \multicolumn{1}{ |c|}{[3]} & [9]\\
       \hline 
\multicolumn{1}{|c|}{-14}  & \multicolumn{1}{ c|}{[3]} & [3] & \multicolumn{1}{ |c|}{[]}  & \multicolumn{1}{ c|}{[9]} & \multicolumn{1}{ c|}{[3]} & [3] &[27] & \multicolumn{1}{ |c|}{[3]} & [9]\\
       \hline 
\multicolumn{1}{|c|}{67}  & \multicolumn{1}{ c|}{[3]} & [3] & \multicolumn{1}{ |c|}{[]}  & \multicolumn{1}{ c|}{[3]} & \multicolumn{1}{ c|}{[]} & [3] &[3] & \multicolumn{1}{ |c|}{[]} & [3]\\
       \hline 
\end{tabular}
\end{center}

In the next table we present all the cases meeting with the following criteria: 

$-100<d<3\,000$, 3 splits in $\Q(\sqrt{d})$ and the groups $\Clog{K}$ and $C\ell'_{K}$ are not trivial.

\begin{center}
\begin{tabular}{c|ccc|ccc|ccc|}
\cline{2-10}
  & \multicolumn{3}{c|}{\multirow{2}{*}{$\Q(\sqrt{d})$}} & \multicolumn{3}{c|}{\multirow{2}{*}{$\Q(\sqrt{d},\cos(2\pi/9))$}} & \multicolumn{3}{c|}{\multirow{2}{*}{$\Q(\sqrt{d},\cos(2\pi/27))$}}\\
   & & & & & & & & & \\
  \hline
   \multicolumn{1}{ |c|}{\multirow{2}{*}{$d$}}     & \multicolumn{1}{ c|}{\multirow{2}{*}{$\Clog{K}$}} & \multirow{2}{*}{$\Clog{K}^{[\ell]}$} & \multicolumn{1}{ |c|}{\multirow{2}{*}{$C\ell'_{K}$}} & \multicolumn{1}{c|}{\multirow{2}{*}{$\Clog{K_{1}}$}} & \multicolumn{1}{ c|}{\multirow{2}{*}{$\Clog{K_{1}}^{[\ell]}$}} & \multirow{2}{*}{$C\ell'_{K_{1}}$} & \multirow{2}{*}{$\Clog{K_{2}}$} & \multicolumn{1}{ |c|}{\multirow{2}{*}{$\Clog{K_{2}}^{[\ell]}$}} & \multirow{2}{*}{$C\ell'_{K_{2}}$}  \\
   
    \multicolumn{1}{ |c|}{} & \multicolumn{1}{ c|}{} &  & \multicolumn{1}{ |c|}{} & \multicolumn{1}{ c|}{}  & \multicolumn{1}{ c|}{} & \multicolumn{1}{ c|}{} & & \multicolumn{1}{ |c|}{} & \\
  \hline
  \multicolumn{1}{|c|}{1129}  & \multicolumn{1}{ c|}{[3]} & [] & \multicolumn{1}{ |c|}{[3]}  & \multicolumn{1}{ c|}{[3]} & \multicolumn{1}{ c|}{[]} & [3] &[3] & \multicolumn{1}{ |c|}{[]} & [3]\\
       \hline 
\multicolumn{1}{|c|}{1654}  & \multicolumn{1}{ c|}{[3]} & [] & \multicolumn{1}{ |c|}{[3]}  & \multicolumn{1}{ c|}{[3]} & \multicolumn{1}{ c|}{[]} & [3] &[3] & \multicolumn{1}{ |c|}{[]} & [3]\\
       \hline 
\multicolumn{1}{|c|}{1954}  & \multicolumn{1}{ c|}{[3]} & [] & \multicolumn{1}{ |c|}{[3]}  & \multicolumn{1}{ c|}{[3]} & \multicolumn{1}{ c|}{[]} & [3] &[3] & \multicolumn{1}{ |c|}{[]} & [3]\\
       \hline 
\multicolumn{1}{|c|}{2419}  & \multicolumn{1}{ c|}{[3]} & [] & \multicolumn{1}{ |c|}{[3]}  & \multicolumn{1}{ c|}{[3]} & \multicolumn{1}{ c|}{[]} & [3] &[3] & \multicolumn{1}{ |c|}{[]} & [3]\\
       \hline 
\multicolumn{1}{|c|}{2713}  & \multicolumn{1}{ c|}{[3]} & [] & \multicolumn{1}{ |c|}{[3]}  & \multicolumn{1}{ c|}{[9]} & \multicolumn{1}{ c|}{[]} & [9] &[9] & \multicolumn{1}{ |c|}{[]} & [9]\\
       \hline 
\multicolumn{1}{|c|}{2971}  & \multicolumn{1}{ c|}{[3]} & [] & \multicolumn{1}{ |c|}{[3]}  & \multicolumn{1}{ c|}{[3]} & \multicolumn{1}{ c|}{[]} & [3] &[3] & \multicolumn{1}{ |c|}{[]} & [3]\\
       \hline 
\end{tabular}
\end{center}

In the next table we present all the cases issued with the following criteria: 

$-100<d<20\,000$, 3 splits in $\Q(\sqrt{d})$ and the groups $\Clog{K}^{[\ell]}$ and $C\ell'_{K}$ are non trivial.

\begin{center}
\begin{tabular}{c|ccc|ccc|ccc|}
\cline{2-10}
  & \multicolumn{3}{c|}{\multirow{2}{*}{$\Q(\sqrt{d})$}} & \multicolumn{3}{c|}{\multirow{2}{*}{$\Q(\sqrt{d},\cos(2\pi/9))$}} & \multicolumn{3}{c|}{\multirow{2}{*}{$\Q(\sqrt{d},\cos(2\pi/27))$}}\\
   & & & & & & & & & \\
  \hline
   \multicolumn{1}{ |c|}{\multirow{2}{*}{$d$}}     & \multirow{2}{*}{$\Clog{K}$} & \multicolumn{1}{ |c}{\multirow{2}{*}{$\Clog{K}^{[\ell]}$}} & \multicolumn{1}{ |c|}{\multirow{2}{*}{$C\ell'_{K}$}} & \multicolumn{1}{ c|}{\multirow{2}{*}{$\Clog{K_{1}}$}} & \multicolumn{1}{ c|}{\multirow{2}{*}{$\Clog{K_{1}}^{[\ell]}$}} & \multirow{2}{*}{$C\ell'_{K_{1}}$} & \multirow{2}{*}{$\Clog{K_{2}}$} & \multicolumn{1}{ |c|}{\multirow{2}{*}{$\Clog{K_{2}}^{[\ell]}$}} & \multirow{2}{*}{$C\ell'_{K_{2}}$}  \\
   
    \multicolumn{1}{ |c|}{} & & \multicolumn{1}{ |c}{} & \multicolumn{1}{ |c|}{} & \multicolumn{1}{ c|}{}  & \multicolumn{1}{ c|}{} & \multicolumn{1}{ c|}{} & & \multicolumn{1}{ |c|}{} & \\
  \hline
  \multicolumn{1}{|c|}{3739}  & \multicolumn{1}{ c|}{[9]} & [3] & \multicolumn{1}{ |c|}{[3]}  & \multicolumn{1}{ c|}{[27]} & \multicolumn{1}{ c|}{[3]} & [9] &[81] & \multicolumn{1}{ |c|}{[3]} & [27]\\
       \hline 
\multicolumn{1}{|c|}{7726}  & \multicolumn{1}{ c|}{[3,3]} & [3] & \multicolumn{1}{ |c|}{[3]}  & \multicolumn{1}{ c|}{[9,3,3]} & \multicolumn{1}{ c|}{[]} & [9,3,3] &[27,3,3] & \multicolumn{1}{ |c|}{[]} & [27,3,3]\\
       \hline 
\multicolumn{1}{|c|}{11545}  & \multicolumn{1}{ c|}{[9]} & [3] & \multicolumn{1}{ |c|}{[3]}  & \multicolumn{1}{ c|}{[9]} & \multicolumn{1}{ c|}{[]} & [9] &[9] & \multicolumn{1}{ |c|}{[]} & [9]\\
       \hline 
\multicolumn{1}{|c|}{17134}  & \multicolumn{1}{ c|}{[9]} & [3] & \multicolumn{1}{ |c|}{[3]}  & \multicolumn{1}{ c|}{[9]} & \multicolumn{1}{ c|}{[]} & [9] &[9] & \multicolumn{1}{ |c|}{[]} & [9]\\
       \hline 
\multicolumn{1}{|c|}{19330}  & \multicolumn{1}{ c|}{[9]} & [3] & \multicolumn{1}{ |c|}{[3]}  & \multicolumn{1}{ c|}{[9]} & \multicolumn{1}{ c|}{[]} & [9] &[9] & \multicolumn{1}{ |c|}{[]} & [9]\\
       \hline 
\end{tabular}
\end{center}
 
\subsection{Logarithmic invariants: Cyclotomic case}
 
We are now interested in explicitly knowing the values of the logarithmic invariants $\mut$ et $\lat$. In \cite{Jaulent86} Jaulent showed (in the case $K_{\infty}/K$ a cyclotomic $\Zl$-extension) that the invariant $\mut$ attached to the logarithmic class group equals its classical counterpart $\mu$. Hence we can state the following result. 
 
\begin{theorem}[Logarithmic Ferrero-Washington] Let $K$ be an abelian number field over $\Q$, and let $K^{c}/K$ be its cyclotomic $\Zl$-extension. Then
$$\mut=0.$$
\end{theorem}

\begin{remark} 
We have $\mut=0$ for the examples given in the preceding section \ref{sub:quadfiels}. 
\end{remark}

Also in \cite[IV.1]{Jaulent86} it is shown that $\lambda'=\lat$, where $\lambda'$ is the invariant attached to the $\ell$-subgroup of $\ell$-classes $C\ell'_{n}$, indeed this follows from the fact that in the cyclotomic tower we have $\Clog{n}\isom \Cc'/\w_{n}\Cc'$. Furthermore Gold's papers \cite{Gold74Exa1} and \cite{Gold74Exa2} allow us to compute explicitly the value of $\lat$ for imaginary quadratic fields:

\begin{theorem}[\cite{Gold74Exa1,Gold74Exa2}]\label{thm:Gold2}
Let $K$ be an imaginary quadratic number field and let $K_{n}$ be the $n$-th layer of the cyclotomic $\Zl$-extension of $K$. Let $e'_{n}$ be the power of the order of the $\ell$-subgroup of $\ell$-classes $C\ell'_{n}$ attached to $K_{n}$. If there is $n\geq 1$, such that $e'_{n}-e'_{n-1}<\f(\ell^{n})$, then $$\lambda'=e'_{n}-e'_{n-1}.$$
\end{theorem}

We use these facts to compute $\mut$ and $\lat$ for the layers $K$, $K_{1}$ and $K_{2}$ of the cyclotomic $\Z_{3}$-extension of the imaginary quadratic number field $\Q(\sqrt{d})$ for $-100\leq d \leq -1$. 

\begin{center}
\begin{tabular}{c|cc|cc|cc|c|c|cc|}
\cline{2-11}
  & \multicolumn{2}{c|}{\multirow{2}{*}{$\Q(\sqrt{d})$}} & \multicolumn{2}{c|}{\multirow{2}{*}{$\Q(\sqrt{d},\cos(2\pi/9))$}} & \multicolumn{2}{c|}{\multirow{2}{*}{$\Q(\sqrt{d},\cos(2\pi/27))$}} & \multirow{2}{*}{$\f(3)=2$} & \multirow{2}{*}{$\f(9)=6$} & \multicolumn{2}{c|}{\multirow{2}{*}{$\mut=0$}}\\
   & & & & & & & & & & \\
  \hline
   \multicolumn{1}{ |c|}{\multirow{2}{*}{$d$}}     & \multirow{2}{*}{$\Clog{K}$}  & \multicolumn{1}{ |c|}{\multirow{2}{*}{$C\ell'_{K}$}} & \multirow{2}{*}{$\Clog{K_{1}}$}  & \multicolumn{1}{ |c|}{\multirow{2}{*}{$C\ell'_{K_{1}}$}} & \multirow{2}{*}{$\Clog{K_{2}}$} & \multicolumn{1}{ |c|}{\multirow{2}{*}{$C\ell'_{K_{2}}$}} & \multirow{2}{*}{$e'_{1}-e'_{0}$} & \multicolumn{1}{ c|}{\multirow{2}{*}{$e'_{2}-e'_{1}$}} & \multirow{2}{*}{$\lat$} & \multicolumn{1}{ |c|}{\multirow{2}{*}{$\widetilde{\nu}$}} \\
   
    \multicolumn{1}{ |c|}{} & & \multicolumn{1}{ |c|}{} & & \multicolumn{1}{ |c|}{}  & \multicolumn{1}{ c|}{} & \multicolumn{1}{ c|}{} & & \multicolumn{1}{ c|}{} & & \multicolumn{1}{ |c|}{} \\
  \hline
 \multicolumn{1}{ |c|}{-87} & [3] & \multicolumn{1}{ |c|}{[3]} & [9] & \multicolumn{1}{ |c|}{[9]}  & \multicolumn{1}{ c|}{[27]} & \multicolumn{1}{ c|}{[27]} & 1 & \multicolumn{1}{ c|}{1} & 1 & \multicolumn{1}{ |c|}{1} \\
       \hline 
\multicolumn{1}{ |c|}{-86} & [3] & \multicolumn{1}{ |c|}{[]} & [9, 9] & \multicolumn{1}{ |c|}{[9, 3]}  & \multicolumn{1}{ c|}{[27, 27]} & \multicolumn{1}{ c|}{[27, 9]} & 3 & \multicolumn{1}{ c|}{2} & 2 & \multicolumn{1}{ |c|}{2} \\
       \hline 
\multicolumn{1}{ |c|}{-74} & [9] & \multicolumn{1}{ |c|}{[]} & [27] & \multicolumn{1}{ |c|}{[3]}  & \multicolumn{1}{ c|}{[81]} & \multicolumn{1}{ c|}{[9]} & 1 & \multicolumn{1}{ c|}{1} & 1 & \multicolumn{1}{ |c|}{2} \\
       \hline 
\multicolumn{1}{ |c|}{-65} & [3] & \multicolumn{1}{ |c|}{[]} & [9] & \multicolumn{1}{ |c|}{[3]}  & \multicolumn{1}{ c|}{[27]} & \multicolumn{1}{ c|}{[9]} & 1 & \multicolumn{1}{ c|}{1} & 1 & \multicolumn{1}{ |c|}{1} \\
       \hline 
\multicolumn{1}{ |c|}{-61} & [3] & \multicolumn{1}{ |c|}{[3]} & [9] & \multicolumn{1}{ |c|}{[9]}  & \multicolumn{1}{ c|}{[27]} & \multicolumn{1}{ c|}{[27]} & 1 & \multicolumn{1}{ c|}{1} & 1 & \multicolumn{1}{ |c|}{1} \\
       \hline 
\multicolumn{1}{ |c|}{-47} & [9] & \multicolumn{1}{ |c|}{[]} & [27] & \multicolumn{1}{ |c|}{[3]}  & \multicolumn{1}{ c|}{[81]} & \multicolumn{1}{ c|}{[9]} & 1 & \multicolumn{1}{ c|}{1} & 1 & \multicolumn{1}{ |c|}{2} \\
       \hline 
\multicolumn{1}{ |c|}{-41} & [27] & \multicolumn{1}{ |c|}{[]} & [81, 3] & \multicolumn{1}{ |c|}{[3, 3]}  & \multicolumn{1}{ c|}{[243, 9]} & \multicolumn{1}{ c|}{[9, 9]} & 2 & \multicolumn{1}{ c|}{2} & 2 & \multicolumn{1}{ |c|}{3} \\
       \hline 
\multicolumn{1}{ |c|}{-35} & [3] & \multicolumn{1}{ |c|}{[]} & [9] & \multicolumn{1}{ |c|}{[3]}  & \multicolumn{1}{ c|}{[27]} & \multicolumn{1}{ c|}{[9]} & 1 & \multicolumn{1}{ c|}{1} & 1 & \multicolumn{1}{ |c|}{1} \\
       \hline 
\multicolumn{1}{ |c|}{-31} & [3] & \multicolumn{1}{ |c|}{[3]} & [9] & \multicolumn{1}{ |c|}{[9]}  & \multicolumn{1}{ c|}{[27]} & \multicolumn{1}{ c|}{[27]} & 1 & \multicolumn{1}{ c|}{1} & 1 & \multicolumn{1}{ |c|}{1} \\
       \hline 
\multicolumn{1}{ |c|}{-14} & [3] & \multicolumn{1}{ |c|}{[]} & [9] & \multicolumn{1}{ |c|}{[3]}  & \multicolumn{1}{ c|}{[27]} & \multicolumn{1}{ c|}{[9]} & 1 & \multicolumn{1}{ c|}{1} & 1 & \multicolumn{1}{ |c|}{1} \\
       \hline 
\end{tabular}
\end{center}

We now compute the logarithmic class groups for the layers $K$ and $K_{1}$ of the cyclotomic $\Z_{5}$ and $\Z_{7}$ extensions of the imaginary quadratic number fields $\Q(\sqrt{d})$ for $-100\leq d \leq -1$. In all these cases, it suffices to compute $C\ell_{K}$ and $C\ell_{K_{1}}$ to know the logarithmic invariants. 

\textit{Case $\ell=5$}

\begin{center}
\begin{tabular}{c|cc|cc|c|cc|}
\cline{2-8}
  & \multicolumn{2}{c|}{\multirow{2}{*}{$K=\Q(\sqrt{d})$}} & \multicolumn{2}{c|}{\multirow{2}{*}{$K_{1}$}}  & \multirow{2}{*}{$\f(5)=4$}  & \multicolumn{2}{c|}{\multirow{2}{*}{$\mut=0$}}\\
   & & & & & & &  \\
  \hline
  \multicolumn{1}{ |c|}{\multirow{2}{*}{$d$}}     & \multirow{2}{*}{$\Clog{K}$}  & \multicolumn{1}{ |c|}{\multirow{2}{*}{$C\ell'_{K}$}} & \multirow{2}{*}{$\Clog{K_{1}}$}  & \multicolumn{1}{ |c|}{\multirow{2}{*}{$C\ell'_{K_{1}}$}}  & \multirow{2}{*}{$e'_{1}-e'_{0}$}  & \multirow{2}{*}{$\lat$} & \multicolumn{1}{ |c|}{\multirow{2}{*}{$\widetilde{\nu}$}} \\
   
    \multicolumn{1}{ |c|}{} & & \multicolumn{1}{ |c|}{} & & \multicolumn{1}{ |c|}{}  & \multicolumn{1}{ c|}{} &  & \multicolumn{1}{ |c|}{} \\
  \hline
  \multicolumn{1}{ |c|}{-51} & [125] & \multicolumn{1}{ |c|}{[]} & [625] & \multicolumn{1}{ |c|}{[5]}  &  1 &  \multicolumn{1}{ c|}{1} & 3 \\
       \hline 
\multicolumn{1}{ |c|}{-47} & [5] & \multicolumn{1}{ |c|}{[5]} & [25] & \multicolumn{1}{ |c|}{[25]}  &  1 &  \multicolumn{1}{ c|}{1} & 1 \\
       \hline 
\multicolumn{1}{ |c|}{-41} & [5] & \multicolumn{1}{ |c|}{[]} & [25] & \multicolumn{1}{ |c|}{[5]}  &  1 &  \multicolumn{1}{ c|}{1} & 1 \\
       \hline 
\multicolumn{1}{ |c|}{-34} & [5] & \multicolumn{1}{ |c|}{[]} & [25] & \multicolumn{1}{ |c|}{[5]}  &  1 &  \multicolumn{1}{ c|}{1} & 1 \\
       \hline 
\multicolumn{1}{ |c|}{-26} & [5] & \multicolumn{1}{ |c|}{[]} & [25] & \multicolumn{1}{ |c|}{[5]}  &  1 &  \multicolumn{1}{ c|}{1} & 1 \\
       \hline 
\multicolumn{1}{ |c|}{-11} & [5] & \multicolumn{1}{ |c|}{[]} & [25] & \multicolumn{1}{ |c|}{[5]}  &  1 &  \multicolumn{1}{ c|}{1} & 1 \\
       \hline 
\end{tabular}
\end{center}

\textit{Case $\ell=7$}\\

\begin{center}
\begin{tabular}{c|cc|cc|c|cc|}
\cline{2-8}
  & \multicolumn{2}{c|}{\multirow{2}{*}{$K=\Q(\sqrt{d})$}} & \multicolumn{2}{c|}{\multirow{2}{*}{$K_{1}$}}  & \multirow{2}{*}{$\f(7)=6$}  & \multicolumn{2}{c|}{\multirow{2}{*}{$\mut=0$}}\\
   & & & & & & &  \\
  \hline
  \multicolumn{1}{ |c|}{\multirow{2}{*}{$d$}}     & \multirow{2}{*}{$\Clog{K}$}  & \multicolumn{1}{ |c|}{\multirow{2}{*}{$C\ell'_{K}$}} & \multirow{2}{*}{$\Clog{K_{1}}$}  & \multicolumn{1}{ |c|}{\multirow{2}{*}{$C\ell'_{K_{1}}$}}  & \multirow{2}{*}{$e'_{1}-e'_{0}$}  & \multirow{2}{*}{$\lat$} & \multicolumn{1}{ |c|}{\multirow{2}{*}{$\widetilde{\nu}$}} \\
   
    \multicolumn{1}{ |c|}{} & & \multicolumn{1}{ |c|}{} & & \multicolumn{1}{ |c|}{}  & \multicolumn{1}{ c|}{} &  & \multicolumn{1}{ |c|}{} \\
  \hline
\multicolumn{1}{ |c|}{-89} & [7] & \multicolumn{1}{ |c|}{[]} & [49] & \multicolumn{1}{ |c|}{[7]}  &  1 &  \multicolumn{1}{ c|}{1} & 1 \\
       \hline 
\multicolumn{1}{ |c|}{-73} & [7] & \multicolumn{1}{ |c|}{[]} & [49] & \multicolumn{1}{ |c|}{[7]}  &  1 &  \multicolumn{1}{ c|}{1} & 1 \\
       \hline 
\multicolumn{1}{ |c|}{-71} & [7] & \multicolumn{1}{ |c|}{[7]} & [49] & \multicolumn{1}{ |c|}{[49]}  &  1 &  \multicolumn{1}{ c|}{1} & 1 \\
       \hline 
\multicolumn{1}{ |c|}{-34} & [7] & \multicolumn{1}{ |c|}{[]} & [49] & \multicolumn{1}{ |c|}{[7]}  &  1 &  \multicolumn{1}{ c|}{1} & 1 \\
       \hline 
\end{tabular}
\end{center}

\subsection{Logarithmic invariants: Non-cyclotomic case}

We now compute logarithmic invariants of some anti-cyclotomic extensions. For that we extrapolate the computations made by Hubbard and Washington \cite{Hubbard&Washington17}.

Let $K=\Q(\sqrt{-3})$ and $\ell=3$. Let $K_{\infty}$ be the anti-cyclotomic $\Zl$-extension of $K$, that is $\Gal(K/\Q)$ acts by -1 on $\Gal(K_{\infty}/K)$. Let $K_{1}$ be the first layer of $K_{\infty}$, i.e. $K_{1}=K(\sqrt[3]{3})$. The extension $L=\Q(\sqrt{-3},\sqrt[3]{d})$ is an extension of degree 3 of $K$. Let $L_{n}$ be such that $L_{n}=LK_{n}$, therefore $L_{\infty}/L$ is a non-cyclotomic $\Zl$-extension.

All places of $L$ above 3 are ramified in $L_{\infty}/L$, and hence they split finitely in the sense of Definition \ref{def:DecompoFinie}. Hence Theorem \ref{Thm:relinva} is valid for our examples, this implies that $\tilde{\mu}=\mu$. For the examples below we know the values of $\mu$, $\lambda$ and $\nu$ \cite[\S 9.1]{Hubbard&Washington17}. Furthermore we know that for $n$ big enough we have
\[\tilde{e}_{n} =  \tilde{\mu} \ell^{n}+ \tilde{\lambda}n + \tilde{\nu},\]
therefore 
\[\tilde{e}_{n+1}-\tilde{e}_{n}-\tilde{\mu}( \ell^{n+1}-\ell^{n})= \tilde{\lambda}. \]
\textbf{Assumption}: The $n$ in theorem \ref{thm:Iwalog} is the same as in the classical situation (which in the following examples turns out to be $n=0$). 

Then we compute $\tilde{e_{1}}$ and $\tilde{e_{0}}$ with PARI/GP and use the above formula to obtain $\tilde{\lambda}$. 

\begin{remark} It would be interesting to give conditions under which the above assumption is true. In the examples below it seems to work since the ramification indices and the logarithmic ramification indices are equal for all the primes above $3$ in $L$ and $L_{1}$.
\end{remark}

\begin{center}
\begin{tabular}{c|cc|cc|ccccc|}
\cline{2-10}
  & \multicolumn{2}{c|}{\multirow{2}{*}{$L=\Q(\sqrt{-3},\sqrt[3]{d})$}} & \multicolumn{2}{c|}{\multirow{2}{*}{$L_{1}$}} & \multicolumn{5}{c|}{\multirow{2}{*}{Invariants}} \\
   & & & & & & & & &  \\
  \hline
   \multicolumn{1}{ |c|}{\multirow{2}{*}{$d$}}     & \multirow{2}{*}{$\Clog{L}$}  & \multicolumn{1}{ |c|}{\multirow{2}{*}{$C\ell'_{L}$}} & \multirow{2}{*}{$\Clog{L_{1}}$}  & \multicolumn{1}{ |c|}{\multirow{2}{*}{$C\ell'_{L_{1}}$}} & \multirow{2}{*}{$\mu=\mut$} & \multicolumn{1}{ |c|}{\multirow{2}{*}{$\lat$}} & \multirow{2}{*}{$\widetilde{\nu}$} & \multicolumn{1}{| c|}{\multirow{2}{*}{$\lambda$}} & \multirow{2}{*}{$\nu$}  \\
   
    \multicolumn{1}{ |c|}{} & & \multicolumn{1}{ |c|}{} & & \multicolumn{1}{ |c|}{}  & \multicolumn{1}{ c|}{} & \multicolumn{1}{ c|}{} & & \multicolumn{1}{ |c|}{} & \\
  \hline
 
 \multicolumn{1}{ |c|}{10} & [] & \multicolumn{1}{ |c|}{[]} & [3, 3, 3] & \multicolumn{1}{ |c|}{[3, 3]}  & \multicolumn{1}{ c|}{1} & \multicolumn{1}{ c|}{1} & -1 & \multicolumn{1}{ |c|}{0} & -1  \\
       \hline
\multicolumn{1}{ |c|}{22} & [3] & \multicolumn{1}{ |c|}{[3]} & [3, 3, 3] & \multicolumn{1}{ |c|}{[3, 3, 3]}  & \multicolumn{1}{ c|}{1} & \multicolumn{1}{ c|}{0} & 0 & \multicolumn{1}{ |c|}{0} & 1  \\
       \hline
\multicolumn{1}{ |c|}{34} & [3] & \multicolumn{1}{ |c|}{[3]} & [3, 3, 3] & \multicolumn{1}{ |c|}{[3, 3, 3]}  & \multicolumn{1}{ c|}{1} & \multicolumn{1}{ c|}{0} & 0 & \multicolumn{1}{ |c|}{0} & 1  \\
       \hline
\multicolumn{1}{ |c|}{44} & [] & \multicolumn{1}{ |c|}{[]} & [9, 3, 3] & \multicolumn{1}{ |c|}{[3, 3]}  & \multicolumn{1}{ c|}{1} & \multicolumn{1}{ c|}{2} & -1 & \multicolumn{1}{ |c|}{0} & -1  \\
       \hline
\multicolumn{1}{ |c|}{46} & [] & \multicolumn{1}{ |c|}{[]} & [3, 3, 3] & \multicolumn{1}{ |c|}{[3, 3]}  & \multicolumn{1}{ c|}{1} & \multicolumn{1}{ c|}{1} & -1 & \multicolumn{1}{ |c|}{0} & -1  \\
       \hline
\multicolumn{1}{ |c|}{58} & [3] & \multicolumn{1}{ |c|}{[3]} & [3, 3, 3] & \multicolumn{1}{ |c|}{[3, 3, 3]}  & \multicolumn{1}{ c|}{1} & \multicolumn{1}{ c|}{0} & 0 & \multicolumn{1}{ |c|}{0} & 1  \\
       \hline
\multicolumn{1}{ |c|}{68} & [3] & \multicolumn{1}{ |c|}{[3]} & [3, 3, 3] & \multicolumn{1}{ |c|}{[3, 3, 3]}  & \multicolumn{1}{ c|}{1} & \multicolumn{1}{ c|}{0} & 0 & \multicolumn{1}{ |c|}{0} & 1  \\
       \hline
\multicolumn{1}{ |c|}{85} & [3] & \multicolumn{1}{ |c|}{[3]} & [3, 3, 3] & \multicolumn{1}{ |c|}{[3, 3, 3]}  & \multicolumn{1}{ c|}{1} & \multicolumn{1}{ c|}{0} & 0 & \multicolumn{1}{ |c|}{0} & 1  \\
       \hline
\multicolumn{1}{ |c|}{92} & [3] & \multicolumn{1}{ |c|}{[3]} & [3, 3, 3] & \multicolumn{1}{ |c|}{[3, 3, 3]}  & \multicolumn{1}{ c|}{1} & \multicolumn{1}{ c|}{0} & 0 & \multicolumn{1}{ |c|}{0} & 1  \\
       \hline
\multicolumn{1}{ |c|}{110} & [3, 3, 3] & \multicolumn{1}{ |c|}{[3, 3, 3]} & [3, 3, 3, 3, 3, 3, 3] & \multicolumn{1}{ |c|}{[3, 3, 3, 3, 3, 3, 3]}  & \multicolumn{1}{ c|}{2} & \multicolumn{1}{ c|}{0} & 1 & \multicolumn{1}{ |c|}{0} & 2  \\
       \hline
\multicolumn{1}{ |c|}{164} & [3] & \multicolumn{1}{ |c|}{[3]} & [3, 3, 3] & \multicolumn{1}{ |c|}{[3, 3, 3]}  & \multicolumn{1}{ c|}{1} & \multicolumn{1}{ c|}{0} & 0 & \multicolumn{1}{ |c|}{0} & 1  \\
       \hline
\multicolumn{1}{ |c|}{170} & [3] & \multicolumn{1}{ |c|}{[3]} & [3, 3, 3, 3, 3] & \multicolumn{1}{ |c|}{[3, 3, 3, 3, 3]}  & \multicolumn{1}{ c|}{2} & \multicolumn{1}{ c|}{0} & -1 & \multicolumn{1}{ |c|}{0} & 0  \\
       \hline
\multicolumn{1}{ |c|}{230} & [3, 3, 3] & \multicolumn{1}{ |c|}{[3, 3, 3]} & [3, 3, 3, 3, 3, 3, 3] & \multicolumn{1}{ |c|}{[3, 3, 3, 3, 3, 3, 3]}  & \multicolumn{1}{ c|}{2} & \multicolumn{1}{ c|}{0} & 1 & \multicolumn{1}{ |c|}{0} & 2  \\
       \hline
\multicolumn{1}{ |c|}{236} & [3] & \multicolumn{1}{ |c|}{[3]} & [3, 3, 3] & \multicolumn{1}{ |c|}{[3, 3, 3]}  & \multicolumn{1}{ c|}{1} & \multicolumn{1}{ c|}{0} & 0 & \multicolumn{1}{ |c|}{0} & 1  \\
       \hline
\multicolumn{1}{ |c|}{253} & [] & \multicolumn{1}{ |c|}{[]} & [3, 3, 3] & \multicolumn{1}{ |c|}{[3, 3]}  & \multicolumn{1}{ c|}{1} & \multicolumn{1}{ c|}{1} & -1 & \multicolumn{1}{ |c|}{0} & -1  \\
       \hline

\end{tabular}
\end{center}

\bibliographystyle{amsplain}
\bibliography{Biblio}  
\end{document}